\documentclass[conference]{IEEEtran}
\IEEEoverridecommandlockouts
\usepackage{cite}
\usepackage{booktabs}
\usepackage{array}
\usepackage{pdflscape}
\usepackage{amsmath,amssymb,amsfonts}
\usepackage{tikz}
\usetikzlibrary{arrows.meta}
\usepackage{graphicx}
\usepackage{textcomp}
\usepackage{xcolor}
\usepackage{stmaryrd}
\def\BibTeX{{\rm B\kern-.05em{\sc i\kern-.025em b}\kern-.08em
    T\kern-.1667em\lower.7ex\hbox{E}\kern-.125emX}}
\begin{document}

\title{A note on finding optimal cut-offs}

\author{Shravan Mohan\\
\textit{Amazon, Bangalore.}}

\maketitle

\begin{abstract}
This paper addresses the challenge of determining optimal cut-offs for a set of $n$ items with $m$ scores to maximize distinguish-ability. The term \textit{distinguish-ability} is defined as the fraction of item pairs assigned to different buckets, where buckets are determined by the number of cut-offs an item's scores exceed. A brute-force approach to this problem is computationally intractable, with complexity growing exponentially with the number of scores. On the other hand, attempts to solve the problem in the continuous domain lead to local minima, making it unreliable. To overcome these challenges, the problem is formulated as a Integer Quadratic Program (IQP). Since IQPs become computationally difficult with even moderate size problems, a surrogate Integer Linear Program (ILP) is introduced, which can be solved more efficiently for larger instances. In addition to these exact methods, a simple heuristic is proposed that offers a balance between solution quality and computational efficiency. This heuristic iteratively adjusts cut-offs for each score, considering a finite-set of meaningful cut-off points. Computational results provide an empirically evidence for the effectiveness of our proposed methods. 
\end{abstract}

\begin{IEEEkeywords}
Optimal cut-offs, distinguish-ability, Integer Quadratic Programming, Integer Linear Programming
\end{IEEEkeywords}

\section{Introduction}
\noindent Consider the problem of dividing a set of $n$ items, each with $m$ numerical scores into different buckets, based on $m$ cut-offs, one for each score. The natural question then is to define an objective for optimality, and also the related bucketing scheme. Since the aim of cut-offs, typically, is to be able to distinguish between items, a suitable measure of \textit{distinguish-ability} must be a natural choice. The measure of \textit{distinguish-ability} as the ratio of the number item pairs that fall into different buckets, to the total number of pairs. Mathematically, suppose that the matrix of scores is given by $S\in \mathbb{R}^{n\times m}$ and $c$ is the vector of the cut-offs; then, distinguish-ability $d(S, c)$ is defined as:
\begin{equation}
    d(S, c) = \frac{\mbox{\# \textit{of pairs of items which fall in different buckets}}}{\mbox{\# \textit{of total pairs of items}}}.
\end{equation}
Since the total of the number of items across buckets is equal to $n$, we also have that
\begin{equation}
    d(S, c) = \frac{n^2 - \mbox{\textit{sum of squares of the  bucket quantities}}}{\mbox{\# \textit{of total pairs of items}}}.
\end{equation}
And hence, maximizing $d(S, c)$ is equivalent to minimizing the \textit{sum of squares} of the numbers of items across buckets. 
As for the choice of buckets, which is also an important consideration, a natural choice is to bucket items based on the number of cut-offs cleared. Another choice could be based on the exact subset of scores for which the items clear cut-offs. However, the main focus shall be on the former, since the number of buckets is $(m+1)$; in the latter, the number of buckets is $2^m$, which becomes numerically infeasible, in terms of time and space complexity, even for a moderately large cases. 
\\

\noindent cut-off based classification systems can play an important role in many real-world applications where clear, actionable cut-offs are necessary for decision-making and communication. cut-offs create distinct categories that can facilitate standardized responses, resource allocation, or intervention strategies. For instance, in educational assessment, grade boundaries help determine advancement or remediation needs; in healthcare, diagnostic cut-offs guide treatment decisions; and in financial risk assessment (such as loan risk), cut-off points help standardize response protocols. For example, the authors of \cite{biomarkers} present a statistical framework to determine the optimal number and location of cut-off points for continuous biomarkers, using cervical cancer data as an example. Traditionally, biomarkers are split into two categories using a single cut-off, but this study highlights the need for multiple cut-offs to better stratify patients into different risk groups. In an example from the finance domain, the authors of \cite{finance} present a rigourous framework for developing optimal scoring cut-off policies in credit and loan portfolios, particularly when managing multiple, often conflicting business objectives, and when multiple scorecards are available, often comparitively non-dominating, thereby making for a challenging problem. Similarly, the authors of \cite{education} develop a mathematical model to determine optimal cut-offs for personnel selection especially when candidates are evaluated using multiple tests. The novelty lies in optimizing both the individual test cut-offs and the collective decision rule. The authors notice that the cut-off for a test used on its own is not necessarily the same when it's part of a battery of tests. \\

\noindent It is also important to note that this approach differs fundamentally from \textit{k-means} clustering \cite{kmeans} and its explainable variants \cite{explainablekmeans}. While \textit{k-means} algorithm partition items into $k$ groups by minimizing within-cluster distances based on item features or scores, our method specifically focuses on maximizing the number of between-group pairs through the strategic placement of cut-offs. Unlike \textit{k-means}, which allows for flexible cluster boundaries determined by distance metrics, our approach creates groups based on strict score cut-offs, where group membership is determined by how many cut-offs an item's score exceeds. This results in an inherently interpretable grouping mechanism, as items are classified based on clear, discrete criteria rather than relative distances to cluster centroids. In contrast to all the aforementioned references, the present work focuses on maximizing distinguish-ability, and to the best of the author's knowledge, this has not been reported in literature thus far. 
\\

\noindent The computational complexity of this cut-off optimization problem presents interesting challenges. While the problem exhibits characteristics suggesting \textit{NP-hardness}, a formal proof remains elusive despite our attempts at finding a suitable reduction. The \textit{brute search} approach to this problem, while guaranteed to find the optimal solution, faces prohibitive computational complexity. For each of the m scores, the relevant cut-off positions are limited to $n+1$ possibilities: the maximum value, the minimum value, and the $n-1$ points halfway between consecutive values in the sorted sequence. While this is more efficient than considering every possible real number as a potential cut-off, it still leads to $(n+1)^m$ possible combinations that need to be evaluated. This exponential growth makes the approach impractical for real-world applications; for instance, with just 100 items ($n=100$) and 5 scores ($m=5$), one would need to evaluate at least 10 billion combinations. Alternatively, the problem can also be reformulated as a \textit{set function} optimization by observing that any set of cut-offs is effectively dominated by some subset of items. Specifically, for any valid set of cut-offs, one can find a subset of items such that taking the minimum score in each dimension from this subset yields those cut-offs. Thus, instead of searching over all possible cut-off values, one can search over subsets of items, where each subset implicitly defines cut-offs as the minimum scores among its items in each dimension. However, while this reformulation provides another combinatorial optimization perspective on the cut-off selection problem, it is easy to find counterexamples that this set function is \textit{not sub-modular} \cite{submodular}, thereby rendering this formulation infructuous. The lack of submodularity can be shown with a simple example. Consider the 50 items wtih their scores given in Table I. Suppose that
\begin{align}
    A &= \{32, 37, 5, 10, 43, 12, 46, 48, 49, 22, 29\}\nonumber\\
    B &= \{32, 1, 37, 5, 10, 11, 43, 12, 46, 45, 48, 49, 22, 29\}\\
    x &= {24} \nonumber.
\end{align}
With set function 
$$
f(E) = \sum_{i=0}^{m}b_i(c(E)) , ~ E\subseteq S, 
$$
where $c(E)$ is the cut-offs which are dominated by the subset $E$ and $b_i(c(E))$ is the number of items in the $i^{\rm th}$-bin, it can be verified that 
\begin{align}
    342 = f(A\cup x) - f(A) < f(B\cup x) - f(B) = 416.
\end{align}
Similarly, the subsets
\begin{align}
    A &= \{23, 46, 47\}\nonumber\\
    B &= \{2, 3, 5, 7, 10, 14, 17, 18, 20, 23, 24, 25, 26, 27, 29, \\&~~~~~~30, 32, 36, 41, 43, 44, 45, 46, 47, 48, 49\}\nonumber\\
    x &= {34} \nonumber.
\end{align}
prove that the function is not supermodular too; the marginal increments for the subsets above can be verified as
\begin{align}
    472 = f(A\cup x) - f(A) < f(B\cup x) - f(B) = 86.
\end{align}
Alternatively, attempting to solve the problem in the continuous domain, with cut-off positions as optimization variables, leads to a non-convex optimization - this is discussed in the next section. Although this formulation is parsimonious in the number of variables, the landscape prone to local minima, making standard gradient-based methods unreliable. Experimental evidence to this end shall be provided in further sections.
\\

\noindent The problem lends itself to a formulation as a Integer Quadratic Program (IQP) using binary variables. This approach allows us to leverage powerful commercial solvers designed for such problems. However, the computational intensity of solving IQPs grows rapidly with problem size; the reason stated as follows. A common technique to solve such IQPs is the Reformulation-Linearization Technique \cite{rlt}, which transforms the IQP into a Integer Linear Program (ILP). That is, for a Integer Quadratic Program with binary variables $\{x_1, \cdots, x_n\}$, the technique first replaces each quadratic term $x_i x_j$ with a new variable $y_{ij}$. Then, it adds linearization constraints that enforce the relationship between these new variables and the original variables. For binary variables, this includes constraints 
\begin{align}
y_{ij} &\leq x_i \nonumber \\
y_{ij} &\leq x_j \\ 
y_{ij} &\geq x_i + x_j - 1,\nonumber 
\end{align}
completing the conversion to a ILP. While this allows for the use of standard ILP solvers, it comes at the cost of significantly increasing the number of variables. Specifically, as shown above, the number of variables in the resulting ILP grows quadratically with respect to the number of variables in the original IQP. This quadratic growth in problem size can quickly lead to large memory requirements and the solution time also increase dramatically. Consequently, while this approach is viable for moderate-sized problems, it may become impractical for very large instances, necessitating the exploration of alternative solution strategies.
\\

\noindent In pursuit of an efficient solution, an alternative approach using ILP is developed that indirectly addresses the original objective, i.e., minimizes the \textit{range} of the number of items across buckets. As shall be shown later, the range serves as a relaxation to the original IQP. Further, a simple greedy search heurisitic is developed, which iteratively changes the cut-offs looking for the best improvement in distinguish-ability, starting from the median. While these two approaches do not guarantee an optimal solution to the original problem, empirical results consistently demonstrate that it produces near-optimal solutions. The ILP formulation and the greedy search offer significant computational advantages, allowing one to solve much larger problem instances compared to the quadratic formulation. However, establishing a formal approximation bounds for these heuristics also remains an open question. 

\begin{figure}[t]
\centering
\fbox{
\begin{minipage}{0.4\textwidth}
\textit{Suppose $\sigma_j(i)$ is the position of the $j^{\rm th}$ score for the $i^{\rm th}$ item, when those are sorted in ascending order.} 
\begin{align*}
\min~ & ~~~  \sum_{i=1}^{m+1}\boldsymbol{t}^2_i\\
\text{sub to} &\\
&\quad \mathbf{X}_{\sigma_j(n),j} \geq \cdots \geq \mathbf{X}_{\sigma_j(1),j}, ~\forall j\\
&\quad \mathbf{X}_{\sigma_j(k),j} = \mathbf{X}_{\sigma_j(k+1),j}, \forall j, k | \mathbf{S}_{\sigma_j(k),j} = \mathbf{S}_{\sigma_j(k+1),j}\\
&\quad \sum_{k=1}^{m+1} (k-1) \mathbf{Y}_{i,k} = \sum_{k=1}^m \mathbf{X}_{i,k},~\forall i\\
&\quad \sum_{k=1}^{m+1} \mathbf{Y}_{i,k} = 1,~\forall i, ~~ \boldsymbol{1}\times \mathbf{Y}  = \boldsymbol{t}\\
&\quad \mathbf{X}, \mathbf{Y} \in \{0, 1\}^{n \times m}.
\end{align*}
\textit{For the $j^{\rm th}$ score, return the center point of $S_{\sigma_j(i),j}$ and $S_{\sigma_j(i+1),j}$ such that $\boldsymbol{X}_{\sigma_j(i),j} \neq \boldsymbol{X}_{\sigma_j(i+1),j}$. }
\end{minipage}
}
\caption{The Integer Quadratic Programming Problem}
\vspace{-0.5cm}
\label{fig:optimization_IQP}
\end{figure}

\section{The Continuous Domain Setting}
\noindent This section shows the inherent difficulty in modeling the problem in the continuous domain setting. Although there may be a possibility of alternative formulations in the continuous domain, the one presented here arises quite naturally. As before, let $n$ is the number of items and $m$ be the number of scores. Let $c_1, c_2, \cdots, c_m$ be the cut-offs; one for each score. Consider the \textit{logistic} function 
\begin{equation}
\sigma_r(x,c_i) = \frac{1}{1+e^{-r(x-c_i)}},
\end{equation} 
for some $r \gg 1$. This function evaluates to almost 1 for $x > t_i$ and to almost 0 for $x < t_i$, and hence classifying an item as \textit{above cut-off $c_i$} (and evaluating to 1 for $r \gg 1$) if it has a score $s_i > t_i$ and as \textit{below cut-off $c_i$}, otherwise (and evaluating to 0 for $r \gg 1$). The only exception is when $s_i=t_i$, in which case the function would evaluate to 0.5, whatever the value of $r$. One way to get around this exception is to add constraints of the following form:
$$
\sigma_r(s(i,j), t_j)(1-\sigma_r(s(i,j), t_j)) \leq \eta,~\forall i,j,
$$
for $\eta \ll 1$. These constraints ensure that the term $\sigma_r(s(i,j), t_j)$ is within a small range of 0 or 1, and hence never becomes equal to 0.5. However, since the focus of this work is not on this continuous domain formulation of the problem, these constraints will be dropped, from here on, for the sake of brevity. Moving on, the bucket number that an item belongs to can be approximated as
\begin{equation}
c_i = \sum_{j=1}^m \sigma_r(s(i,j), t_j).
\end{equation}
Further with this, the count of items in the $i^{\rm th}$ bucket can be approximated as
\begin{equation}
b_i = \sum_{j=1}^n \mu_i(c_j),
\end{equation}
where $\mu_i(.)$ is defined as 
\begin{equation}
\mu_i(x) = 4\sigma_r(x, i)(1-\sigma_r(x,i)).
\end{equation}
Note that this function evaluates to 1 at $x=i$, and evaluates to almost 0 at all other real numbers (for $r\gg 1$). 
So, the original optimization problem can be approximated as:
\begin{equation}
\min_{c_1, c_2, \dots, c_m} \sum_{i=1}^{m+1}b_i^2.
\end{equation}
Although the derivation of the cost in the optimization problem was paved through approximations, the approximation can be \textit{arbitrarily tight} for a given set of cut-offs (\textit{w.r.t.} to the original optimization cost) by increasing $r$.\\

\noindent The optimization routine employed here uses gradient descent with an adaptive learning rate. At each iteration, the gradient is computed (using autograd \cite{autograd}) at the current point $x$. The algorithm initializes the learning rate ($\epsilon$) at 1 and then systematically reduces it by a factor of 1.1 until
$$
d(S, x - \epsilon \nabla d(S, x)) \geq d(S, x).
$$
This adaptive step ensures a reduction in the objective function at each iteration while preventing overshooting. Once an appropriate learning rate is found, the algorithm updates $x$, i.e., 
$$
x \rightarrow x - \epsilon \nabla d(S, x).
$$
This process continues until a critical point is reached, signifying either a local minimum or a point where the gradient becomes negligibly small. As shall be shown later, the solution obtained using this approach does not guarantee optimality. In fact, the local optima obtained this way also does not compare favourably with the other techniques. 
\begin{figure}[t]
\centering
\fbox{
\begin{minipage}{0.45\textwidth}
\textit{Suppose $\sigma_j(i)$ is the position of the $j^{\rm th}$ score for the $i^{\rm th}$ item, when those are sorted in ascending order.}
\begin{align*}
\min~ & ~~~s - t \\
\text{sub to} &\\
&\quad \mathbf{X}_{\sigma_j(n),j} \geq \cdots \geq \mathbf{X}_{\sigma_j(1),j}, ~\forall j\\
&\quad \mathbf{X}_{\sigma_j(k),j} = \mathbf{X}_{\sigma_j(k+1),j}, \forall j, k | \mathbf{S}_{\sigma_j(k),j} = \mathbf{S}_{\sigma_j(k+1),j}\\
&\quad \sum_{k=1}^{m+1} (k-1) \mathbf{Y}_{i,k} = \sum_{k=1}^m \mathbf{X}_{i,k},~\forall i\\
&\quad \sum_{k=1}^{m+1} \mathbf{Y}_{i,k} = 1,~\forall i, ~~ t\leq \mathbf{Y}\times \boldsymbol{1} \leq s\\
&\quad \mathbf{X}, \mathbf{Y} \in \{0, 1\}^{n \times m}.
\end{align*}
\textit{For the $j^{\rm th}$ score, return the center point of $S_{\sigma_j(i),j}$ and $S_{\sigma_j(i+1),j}$ such that $\boldsymbol{X}_{\sigma_j(i),j} \neq \boldsymbol{X}_{\sigma_j(i+1),j}$. }
\end{minipage}
}
\caption{The Integer Linear Programming Heuristic}
\vspace{-0.5cm}
\label{fig:optimization_ILP}
\end{figure}
\section{The Integer Quadratic Programming Approach}
\noindent The formulation of the Integer Quadratic Program to solve the original optimization problem goes as follows. Let one set of the optimization parameters be given by a matrix $\mathbf{X} \in \{0,1\}^{n\times m}$, where 
\begin{equation}
\mathbf{X}_{i,j} = \begin{cases} 1, & \text{if $i^{\rm th}$ item crosses the $j^{\rm th}$ cut-off,} \\ 0, & \text{otherwise.} \end{cases}
\end{equation}
So, it is obvious that if the $i^{\rm th}$ item crosses the $j^{\rm th}$ cut-off, then all items having their respective $j^{\rm th}$ score more than or equal to that of the $i^{\rm th}$ item would also cross the $j^{\rm th}$ cut-off; and this would be true across scores. This implies that
\begin{equation}
\mathbf{X}_{\sigma_j(n)} \geq \mathbf{X}_{\sigma_j(n-1)} \geq \cdots \geq \mathbf{X}_{\sigma_j(2)} \geq \mathbf{X}_{\sigma_j(1)}, ~\forall j,
\end{equation}
where the mapping $\sigma:\llbracket1,n\rrbracket\rightarrow \llbracket1,n\rrbracket$ is such that $\left[S_{\sigma_j(1)}, \cdots, S_{\sigma_j(n)}, \right]$ is the ascending sorted sequence of the $j^{\rm th}$ column of $S$. And since for any score, two items with the same value should either both clear cut-off or both should not clear cut-off, the following set of additional set of constraints also becomes necessary:
\begin{equation}
\quad \mathbf{X}_{\sigma_j(k),j} = \mathbf{X}_{\sigma_j(k+1),j}, \forall j, k | \mathbf{S}_{\sigma_j(k),j} = \mathbf{S}_{\sigma_j(k+1),j}.
\end{equation}

\noindent Now, consider another set of variables given by a matrix $\textbf{Y} \in \{0,1\}^{n\times (m+1)}$, where
\begin{equation}
\mathbf{Y}_{i,j} = \begin{cases} 1, & \text{if $i^{\rm th}$ item crosses the exactly $j$ cut-offs,} \\ 0, & \text{otherwise.} \end{cases}
\end{equation}
It is immediately clear that
\begin{equation}
    \sum_{k=1}^{m+1}(k-1) \mathbf{Y}_{i,k} = \sum_{k=1}^m \mathbf{X}_{i,k},~\forall i,
\end{equation}
and that 
\begin{equation}
    \sum_{k=1}^{m+1} \mathbf{Y}_{i,k} = 1,~\forall i.
\end{equation}
The column sums of $\boldsymbol{Y}$ give the number of items across buckets. With additional variables $\boldsymbol{t} = [t_1, \cdots, t_n]$ representing the column sums, and the resulting constraint that
\begin{equation}
\boldsymbol{1}\times \boldsymbol{Y} = \boldsymbol{t},
\end{equation}
the cost function can then be written as
\begin{equation}
\sum_{i=1}^{m+1}\boldsymbol{t}^2_i.
\end{equation}
With this, one arrives at the IQP shown in Fig. \ref{fig:optimization_IQP}. \\

\noindent As was mentioned earlier, another way to bucket items can be based on the exact subset of scores for which items clear the cut-offs. Naturally, if there are $m$ scores, there would be $2^m$ possible buckets. It is clear that if one were to make $\boldsymbol{Y}\in \{0,1\}^{n\times 2^m}$, and change the constraints on $Y$ to
\begin{align}
    \sum_{k=1}^{2^m}(k-1) \mathbf{Y}_{i,k} &= \sum_{k=1}^m \mathbf{X}_{i,k},~\forall i,\\
    \sum_{k=1}^{2^m} \mathbf{Y}_{i,k} &= 1,~\forall i,\\
    \boldsymbol{1}\times \boldsymbol{Y} &= \boldsymbol{t},
\end{align}
and the cost function to
\begin{equation}
    \sum_{i=1}^{2^m}\boldsymbol{t}^2_i,
\end{equation}
the resulting IQP would be the formulation to solve the original problem with $2^m$ buckets. However, since the number of buckets increases exponentially with $m$, the formulation can only serve cases with small number of scores ($m \leq 5$). 

\begin{figure}[t]
\centering
\fbox{
\begin{minipage}{0.45\textwidth}
\begin{itemize}
    \item Set each cut-off to \textbf{\textit{median}} of its score.

    \item Calculate \textbf{\textit{initial}} distinguish-ability.
    
    \item For \textbf{\textit{each}} score $i$ (1 to $m$):
    \begin{itemize}
        \item For \textbf{\textit{each}} possible cut-off value $v$:
        \begin{itemize}
            \item Set cut-off to $v$.
            \item Calculate \textbf{\textit{new}} distinguish-ability.
            \item If \textbf{\textit{better}}:
            \begin{itemize}
                \item Keep \textbf{\textit{new}} cut-off.
                \item Update \textbf{\textit{best}} distinguish-ability.
                \item Mark \textbf{\textit{improvement found}}.
            \end{itemize}
        \end{itemize}
    \end{itemize}
    \item If \textbf{\textit{no improvement}}, exit loop. Else \textbf{\textit{continue}}.  
    \item Return \textbf{\textit{final cut-offs}}.
\end{itemize}
\end{minipage}
}
\caption{A Greedy Search Heuristic}
\vspace{-0.5cm}
\label{fig:optimization_greedy}
\end{figure}

\section{{The Integer Linear Programming Approach}}
\noindent To start with the ILP formulation, consider the following lemma.
\\
\noindent \textbf{Lemma}: For $\boldsymbol{x} \in \mathbb{R}^n$, under the constraint that 
\begin{equation}
\sum_{i=1}^n x_i = M,
\label{eqn:equal_sum}
\end{equation}
the range function given by  
$$
\max(\boldsymbol{x}) - \min(\boldsymbol{x})
$$ 
serves a relaxation for 
$$
\sum_{i=1}^{n}x_i^2.
$$.
\\
\noindent \textbf{Proof}: Note that 
$$
\sum_{i=1}^{n}x_i^2 = n\text{Var}(x) + \frac{M^2}{n}  \leq n\left(\max(\boldsymbol{x}) - \min(\boldsymbol{x})\right)^2 + \frac{M^2}{n}.
$$
Since $M$ and $n$ are constants, the above inequality shows that $\left(\max(\boldsymbol{x}) - \min(\boldsymbol{x})\right)^2$ is a relaxation to the original cost. And, since $(.)^2$ is a monotonic function over $\mathbf{R}^+$, one can use  $\left(\max(\boldsymbol{x})- \min(\boldsymbol{x})\right)$ as the cost instead. $\hfill \square$\\

\noindent Since the total of the number of items across buckets is equal to $n$, irrespective of the cut-offs, the lemma above suggests that minimizing the range of the number of items across buckets can serve as a good surrogate to the cost function used in the IQP formulation given in Fig. \ref{fig:optimization_IQP}. p To that end, replacing the following constraint in Fig. \ref{fig:optimization_IQP} 
\begin{align}
    \boldsymbol{1} \times \mbox{\textbf{Y}} = \boldsymbol{t} \nonumber
\end{align}
with the following
\begin{align}
    t \leq \boldsymbol{1} \times \mbox{\textbf{Y}} \leq s, \nonumber
\end{align}
and the cost function with
\begin{align}
    s - t,
\end{align}
one obtains the ILP shown in Fig. \ref{fig:optimization_ILP}.
\section{The Greedy Search Heuristic}
\noindent This heuristic for optimizing cut-off points begins by initializing the cut-offs to the \textit{median values} of each score, providing a balanced starting point. As mentioned earlier, for each score, the number of meaningful cut-off points is limited to at most $n+1$ values. These include the minimum score, the maximum score, and the midpoints between consecutive values in the sorted list of scores.\\

\noindent The method then proceeds \textit{iteratively}, focusing on one score at a time while keeping the others fixed. For each score, it explores this limited set of possible cut-off values, evaluating the distinguish-ability metric for each. After examining all scores, the method selects the modification that yields the \textit{highest increase} in distinguish-ability and updates the cut-off accordingly. This process is repeated in subsequent iterations, continuously refining the cut-offs. The method terminates when an iteration \textit{fails} to produce any improvement in distinguish-ability, indicating that a local optimum has been reached. The method has also been outlined in Fig. \ref{fig:optimization_greedy}.\\

\noindent This method combines the efficiency of a greedy strategy with the thoroughness of examining all meaningful cut-offs for each score, striking a balance between computational speed and solution quality. While it may not guarantee a global optimum, it often produces high-quality solutions quickly, making it particularly suitable for large-scale applications or scenarios requiring rapid decision-making.

\begin{table}[t]
\centering
\caption{Data Table}
\small
\begin{tabular}{@{}r*{3}{>{\centering\arraybackslash}p{1.5cm}}@{}}
\toprule
Item & Score 1 & Score 2 & Score 3 \\
\midrule
1 & -1.25 & -0.94 & -0.53 \\
2 & 0.72 & 1.10 & 0.08 \\
3 & 0.67 & -0.81 & 0.52 \\
4 & -1.36 & -2.26 & 0.73 \\
5 & 0.77 & -1.51 & -0.57 \\
6 & -0.23 & 0.13 & 1.30 \\
7 & -0.14 & 0.49 & 0.06 \\
8 & 0.90 & -0.07 & -0.27 \\
9 & -1.34 & 1.82 & -0.04 \\
10 & -1.08 & 1.08 & 0.58 \\
11 & -0.16 & 0.09 & -0.86 \\
12 & 0.21 & -0.71 & -1.33 \\
13 & 0.93 & 0.38 & -1.31 \\
14 & 2.05 & 1.06 & -1.30 \\
15 & -0.01 & 0.28 & -0.17 \\
16 & -1.17 & 2.32 & 1.37 \\
17 & 0.38 & 0.45 & 0.46 \\
18 & 0.45 & 0.58 & 1.28 \\
19 & -0.29 & 0.81 & 0.88 \\
20 & 0.46 & -1.14 & 0.78 \\
21 & 0.65 & -0.70 & 1.08 \\
22 & -0.88 & 0.60 & 0.18 \\
23 & 2.03 & -0.16 & 1.61 \\
24 & 0.85 & -0.74 & -1.25 \\
25 & -0.58 & -1.22 & -0.97 \\
26 & 0.58 & 0.92 & 0.01 \\
27 & -0.55 & -0.78 & -0.67 \\
28 & -0.70 & 0.75 & 0.46 \\
29 & 0.75 & 0.11 & 1.78 \\
30 & -0.76 & -0.05 & 1.06 \\
31 & 0.68 & -0.58 & 0.88 \\
32 & 0.15 & -1.05 & 0.23 \\
33 & -1.07 & -0.64 & -0.32 \\
34 & 1.19 & -1.91 & -0.41 \\
35 & 0.34 & -2.33 & -0.26 \\
36 & 0.46 & 1.26 & -0.41 \\
37 & 0.64 & -0.23 & 0.13 \\
38 & 0.67 & 0.30 & -1.26 \\
39 & -0.16 & 0.25 & -0.37 \\
40 & -0.60 & -0.19 & 0.13 \\
41 & 1.22 & 0.97 & -1.19 \\
42 & 0.57 & -0.89 & 1.62 \\
43 & 0.53 & -1.22 & -0.11 \\
44 & 0.17 & -0.14 & -0.07 \\
45 & 1.10 & 0.61 & 0.31 \\
46 & -1.44 & -0.40 & -1.13 \\
47 & 1.91 & 0.00 & -1.57 \\
48 & -1.70 & -1.14 & -0.42 \\
49 & -0.03 & 0.85 & -0.96 \\
50 & -0.20 & 0.22 & -0.04 \\
\bottomrule
\end{tabular}
\label{tab:datatab}
\end{table}

\begin{figure}
    \centering
    \includegraphics[width=3.5in]{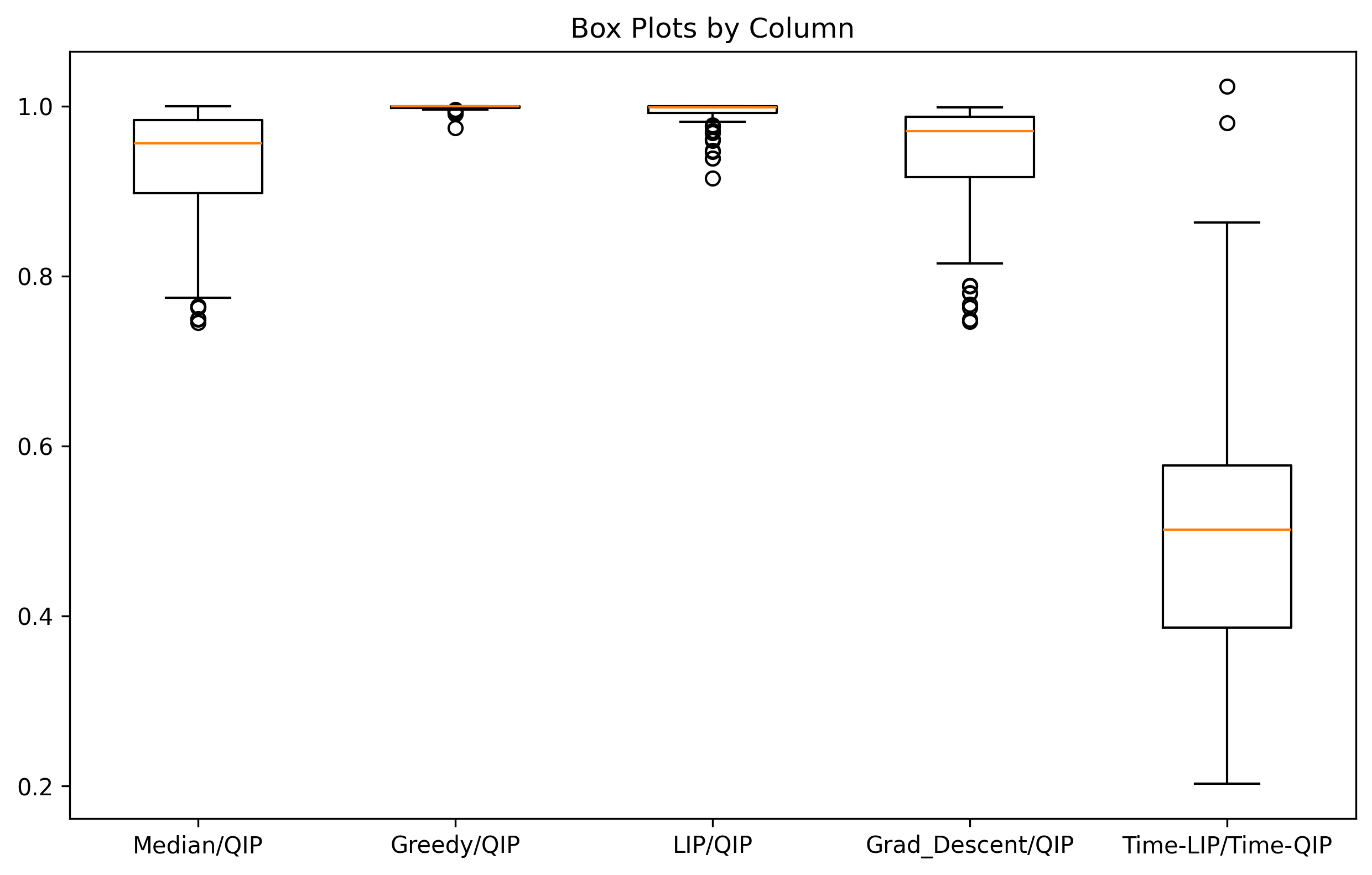}
    \caption{Boxplots of the approximation ratios of the distinguish-abilities obtained for 100 problem instances using (i) median as the cut-offs, (ii) greedy heuristic, (iii) the Integer Linear Programming based relaxation, (iv) gradient descent on the continuous domain approximation. The last boxplot is that of the time-taken to solve the ILP relaxation to the time take to solve the IQP formulation.}
    \label{fig:boxplot}
\end{figure}

\begin{table}
\centering
\caption{cut-off Values, Bucket Sizes, and Distinguishability}
\begin{tabular}{@{}ccccccccc@{}}
\toprule
\multicolumn{3}{c}{cut-off Values} & \multicolumn{4}{c}{\# of Items in bucket} & Distinguishability \\
\cmidrule(lr){1-3} \cmidrule(lr){4-7} \cmidrule(l){8-8}
Score 1 & Score 2 & Score 3 & 0 & 1 & 2 & 3 & Metric \\
\midrule
$-\infty$ &$-\infty$ &$-\infty$ &0&0&0&50&0 \\
$+\infty$ &$+\infty$ &$+\infty$ &50&0&0&0&0 \\
&&&&&&& \\
0.275 &-0.025 & -0.04 & 8 & 14 & 22 & 6 &  0.702 \\
&&&&&&& (\text{Median})\\
&&&&&&& \\
0.36 & 0.415 & -0.015 & 13 & 15 & 17 & 5 & 0.731 \\
&&&&&&& (\text{Greedy})\\
&&&&&&& \\
0.575 & -0.705 & -0.015 & 7 & 17 & 18 & 8 & 0.724 \\
&&&&&&& (\text{ILP})\\
&&&&&&& \\
0.36 & 0.415 & -0.015 & 13 & 15 & 17 & 5 & 0.731 \\
&&&&&&& (\text{IQP})\\
\bottomrule
\end{tabular}
\label{tab:results}
\end{table}

\section{Computational Experiments}
\noindent To begin, consider an example data set with items and their scores is given in Table. \ref{tab:datatab}, and the cut-offs obtained using different techniques shown in Table. \ref{tab:results}. The latter also shows the number of items in the different buckets. It is immediately clear that the median scores are not an effective choice and can be vastly improved upon using the greedy heuristic. In fact, the greedy heuristic ends up giving the optimal solution - the optimal solution is the one obtained by solving the IQP. To evaluate the effectiveness of the proposed approaches, computational experiments across \textit{several} problem instances were solved were solved using multiple methods: the greedy heuristic algorithm, the relaxation-based Integer Linear Program (ILP), the exact Integer Quadratic Program (IQP) formulation. For comparison, solutions (local optima) in the continuous domain were also determined using gradient descent (albeit, \textit{without} any constraints) and a simple median-based cut-off approach. A total of 100 problem instances were chosen; for each problem instance, there were a total of 100 items with 3 scores each. The scores were first generated as a samples from a Gaussian distribution with a randomly generated covariance matrix, and then multiplied by 100 and rounded off to the nearest integer. Integrality ensured that a value of $r=5$ sufficed for the continuous domain setting. Also, the IQP and the ILP were solved using the XPRESS solver \cite{xpress, xpress2}.\\ 

\noindent The results are presented as boxplots in Fig. \ref{fig:boxplot} showing the ratio of distinguish-ability obtained using each method relative to that achieved by the exact IQP formulation. These boxplots reveal that the greedy heuristic \textit{consistently} produces solutions with distinguish-ability comparable to those obtained from the more computationally intensive IQP approach. The ILP-based solutions also performed well (in fact, sometimes better than the greedy heuristic), though with slightly more variability. Additionally, the last boxplot shows the ratio of time taken to solve ILP versus IQP, which demonstrates that the ILP formulation is \textit{substantially faster} than the IQP formulation. In contrast, the continuous domain approach using gradient descent frequently converged to local minima, resulting in suboptimal distinguish-ability, while the median-based cut-off approach performed adequately but was outperformed by the other proposed methods. These results empirically demonstrate that our greedy heuristic offers an excellent balance between computational efficiency and solution quality, making it particularly suitable for large-scale applications.
\section{Conclusion}
\noindent This paper presents a set of approaches to the problem of optimizing cut-off points to maximize distinguish-ability among a set of items with multiple scores. This problem is addressed by first formulating it as a Integer Quadratic Program (IQP) and subsequently as a Integer Linear Program (ILP) with a relaxed cost function. The ILP formulation, in particular, proves effective in solving instances of larger sizes as compared to the IQP approach. The ILP approach was observed to give an approximation ratio of greater than 0.95 on an average, empirically. Further, a heuristic algorithm developed offers a practical alternative for very large instances or scenarios requiring rapid solutions. Computational experiments suggest the effectiveness of this heuristic, with an empirical approximation ratio of greater than 0.95 on an average. Finally, establishing the computational complexity of the problem, and theoretical approximation bounds for the ILP and the heuristic, can be challenging problems for future work.

\end{document}